\documentclass{gtart}

\def\ifplaintex{\expandafter\ifx\csname documentclass\endcsname\relax}

\ifplaintex 
\else
\fi

\def\gtm{{\mathsurround=0pt\it $\cal G\mskip-2mu$eometry \&\ 
$\cal T\!\!$opology $\cal M\mskip-1mu$onographs}}    

\def\gtp{{\mathsurround=0pt\it $\cal G\mskip-2mu$eometry \&\ 
$\cal T\!\!$opology $\cal P\!$ublications}}  

\def\recd{{\small Received:\qua\receiveddate\ifx\reviseddate\relax
\else\qquad Revised:\qua\reviseddate\fi\par}} 

\def\volumenumber#1{\def\thevolumenumber{#1}}
\def\volumeyear#1{\def\thevolumeyear{#1}}
\def\volumename#1{\def\thevolumename{#1}}
\def\papernumber#1{\def\thepapernumber{#1}}
\def\pagenumbers#1#2{\def\startpage{#1}\def\finishpage{#2}}
\def\published#1{\def\publishdate{#1}}
\def\received#1{\def\receiveddate{#1}}
\def\revised#1{\def\reviseddate{#1}}
\def\accepted#1{\def\accepteddate{#1}}

\long\def\asciiabstract#1{\long\def\theasciiabstract{#1}}

\let\\\par
\let\thevolumenumber\relax\let\thepapernumber\relax
\let\thevolumeyear\relax\let\startpage\relax
\let\finishpage\relax\let\publishdate\relax\let\receiveddate\relax
\let\reviseddate\relax\let\accepteddate\relax\let\theasciititle\relax
\let\theasciiauthors\relax
\let\theasciiabstract\relax

\let\theerratum\relax\let\theasciiemail\relax
\let\theshortauthors\relax\let\theshorttitle\relax

\def\startpage{1}\def\finishpage{15}\def\thepapernumber{77}

\volumenumber{2}
\volumename{Proceedings of the Kirbyfest}
\volumeyear{1999}

\long\def\maketitlep{   

\count0=\startpage

\gtm\nl        
{\small Volume \thevolumenumber: \thevolumename\nl 
\ifx\theerratum\relax\else Erratum \erratumnumber\nl\fi
Pages \startpage--\finishpage\nl}

\vglue 0.1truein   

{\parskip=0pt\leftskip 0pt plus 1fil\def\\{\par\smallskip}{\ifplaintex\large
\else\Large\fi\bf\thetitle}\par\medskip}   
\vglue 0.05truein 

{\parskip=0pt\leftskip 0pt plus 1fil\def\\{\par}{\sc\theauthors}
\par\medskip}%
 
\vglue 0.03truein 

{\small\leftskip 25pt\rightskip 25pt{\bf Abstract}\stdspace\theabstract

{\bf AMS Classification}\stdspace\theprimaryclass
\ifx\thesecondaryclass\relax\else; \thesecondaryclass\fi\par
{\bf Keywords}\stdspace \thekeywords\par}\vglue 7pt

}   

\font\phead=cmsl9 scaled 950
\font=cmsl9 scaled 1050
\font\pnum=cmbx10 scaled 913
\font\lnum=cmbx10 
\font\pfoot=cmsl9 scaled 950
\font=cmsl9 scaled 1050
\ifplaintex
\headline{\vbox to 0pt{\vskip -4.5mm\line{\small\phead\ifnum
\count0=\startpage ISSN 1464-8997 (on line)
1464-8989 (printed) \hfill {\pnum\folio}\else\ifodd\count0\def\\{ }%
\ifx\theshorttitle\relax\thetitle\else\theshorttitle\fi\hfill{\pnum\folio}
\else\def\\{ and }{\pnum\folio}\hfill\ifx\theshortauthors\relax\theauthors
\else\theshortauthors\fi\fi\fi}\vss}}
\footline{\vbox to 0pt{\vglue 0mm\line{\small\pfoot\ifnum\count0=\startpage
Published \publishdate:\qua\copyright\ \gtp\hfill\else
\gtm, Volume \thevolumenumber\ (\thevolumeyear)\hfill\fi}\vss
}}
\else
\makeatletter
\def\@oddhead{{\small\lhead\ifnum\count0=\startpage ISSN 1464-8997 (on line)
1464-8989 (printed) \hfill {\lnum\number\count0}\else\ifodd\count0
\def\\{ }\ifx\theshorttitle\relax \thetitle \else\theshorttitle\fi\hfill
{\lnum\number\count0}\else\def\\{ and }{\lnum\number\count0}
\hfill\ifx\theshortauthors\relax 
\theauthors\else\theshortauthors\fi\fi\fi}}\def\@evenhead{@oddhead}
\def\@oddfoot{\small\lfoot\ifnum\count0=\startpage Published \publishdate:\qua\copyright\ \gtp\hfill\else
\gtm, Volume \thevolumenumber\ (\thevolumeyear)\hfill\fi}
\def\@evenfoot{@oddfoot}
\makeatother
\fi

\let\maketitlepage\maketitlep
\let\makeshorttitle\maketitlepage
\let\maketitle\maketitlepage

\newwrite\gtoutfile
\long\gdef\makeheadfile{  
{\def\\{, }\def\s{ }
\immediate\openout\gtoutfile head.xxx
\immediate\write\gtoutfile{To: math@arxiv.org}
\immediate\write\gtoutfile{Subject: put OR rep NNNNN:ppppp}
\immediate\write\gtoutfile{--text follows this line--}
\immediate\write\gtoutfile{Proxy-for: \ifx\theasciiauthors\relax
\theauthors\else\theasciiauthors\fi\s<\ifx\theasciiemail\relax\theemail\else\theasciiemail\fi>}
\immediate\write\gtoutfile{\noexpand\\}
\immediate\write\gtoutfile{Authors: \ifx\theasciiauthors\relax
\theauthors\else\theasciiauthors\fi}
{\def\\{ }\immediate\write\gtoutfile{Title: \ifx\theasciititle\relax
\thetitle\else\theasciititle\fi}}
\immediate\write\gtoutfile{Subj-class: GT or SG, GR etc}
\immediate\write\gtoutfile{MSC-class: \theprimaryclass\ifx\thesecondaryclass\relax\else, \thesecondaryclass\fi}
\immediate\write\gtoutfile{Journal-ref: Geom. Topol. Monogr. \thevolumenumber\s
(\thevolumeyear) \startpage-\finishpage}
\immediate\write\gtoutfile{Comments: Published by Geometry and Topology Monographs at}
\immediate\write\gtoutfile{\s\s\s  http://www.maths.warwick.ac.uk/gt/GTMon\thevolumenumber/paper\thepapernumber.abs.html}
\immediate\write\gtoutfile{\noexpand\\}
\immediate\write\gtoutfile{}
\ifx\theasciiabstract\relax
\immediate\write\gtoutfile{\theabstract}\else
\immediate\write\gtoutfile{\theasciiabstract}\fi
\immediate\write\gtoutfile{}
\immediate\write\gtoutfile{\noexpand\\}
\immediate\write\gtoutfile{}
\immediate\closeout\gtoutfile}}  

\def\maketitlepage{\maketitlep\makeheadfile}
\let\makeshorttitle\maketitlepage
\let\maketitle\maketitlepage

\volumenumber{4}
\volumename{Invariants of knots and 3-manifolds (Kyoto 2001)}
\volumeyear{2002}
\papernumber{3}
\pagenumbers{29}{41}
\received{30 November 2001}
\revised{4 April 2002}
\accepted{22 July 2002}
\published{19 September 2002}

\usepackage{curvesls,epic}

\newtheorem{thm}{Theorem}[section]
\newtheorem{lem}[thm]{Lemma}
\theoremstyle{definition}
\newtheorem{defn}[thm]{Definition}

\newcommand{\Z}{{\bf Z}}
\newcommand{\half}{{\frac{1}{2}}}

\newcommand{\lplus}
{\raisebox{-13pt}
  {
  \begin{picture}(30,30)
  \put(0,0){\vector(1,1){30}}
  \put(30,0){\line(-1,1){13}}
  \put(13,17){\vector(-1,1){13}}
  \end{picture}
  }
}
\newcommand{\lminus}
{\raisebox{-13pt}
  {
  \begin{picture}(30,30)
  \put(0,0){\line(1,1){13}}
  \put(17,17){\vector(1,1){13}}
  \put(30,0){\vector(-1,1){30}}
  \end{picture}
  }
}
\newcommand{\lzero}
{\raisebox{-13pt}
  {
  \begin{picture}(30,30)
  \put(0,0){\line(1,1){10}}
  \qbezier(10,10)(15,15)(10,20)
  \put(10,20){\vector(-1,1){10}}
  \put(30,0){\line(-1,1){10}}
  \qbezier(20,10)(15,15)(20,20)
  \put(20,20){\vector(1,1){10}}
  \end{picture}
  }
}

\newcommand{\forkY}
{\raisebox{-28pt}
  {
  \begin{picture}(60,60)
  \put(30,30){\bigcircle{60}}
  \put(30,30){\circle*{4}}
  \qbezier(0,30)(10,20)(30,20)
  \qbezier(30,20)(50,20)(60,30)
  \drawline(30,0)(30,20)
  \end{picture}
  }
}
\newcommand{\forkL}
{\raisebox{-28pt}
  {
  \begin{picture}(60,60)
  \put(30,30){\bigcircle{60}}
  \put(30,30){\circle*{4}}
  \drawline(0,30)(30,30)
  \drawline(30,0)(15,30)
  \end{picture}
  }
}
\newcommand{\forkR}
{\raisebox{-28pt}
  {
  \begin{picture}(60,60)
  \put(30,30){\bigcircle{60}}
  \put(30,30){\circle*{4}}
  \drawline(30,30)(60,30)
  \drawline(30,0)(45,30)
  \end{picture}
  }
}

\newcommand{\forkTo}
{\raisebox{-28pt}
  {
  \begin{picture}(60,60)
  \put(30,30){\bigcircle{60}}
  \multiput(15,30)(15,0){3}{\circle*{4}}
  \drawline(15,30)(30,30)
  \drawline(30,0)(23,30)
  \end{picture}
  }
}
\newcommand{\forkoT}
{\raisebox{-28pt}
  {
  \begin{picture}(60,60)
  \put(30,30){\bigcircle{60}}
  \multiput(15,30)(15,0){3}{\circle*{4}}
  \drawline(30,30)(45,30)
  \drawline(30,0)(37,30)
  \end{picture}
  }
}
\newcommand{\forkTplus}
{\raisebox{-28pt}
  {
  \begin{picture}(60,60)
  \put(30,30){\bigcircle{60}}
  \put(0,30){\vector(1,0){40}}
  \put(0,30){\line(1,0){60}}
  \put(30,0){\line(0,1){30}}
  \end{picture}
  }
}
\newcommand{\forkTminus}
{\raisebox{-28pt}
  {
  \begin{picture}(60,60)
  \put(30,30){\bigcircle{60}}
  \put(60,30){\vector(-1,0){40}}
  \put(0,30){\line(1,0){60}}
  \put(30,0){\line(0,1){30}}
  \end{picture}
  }
}
\newcommand{\forkTL}
{\raisebox{-28pt}
  {
  \begin{picture}(60,60)
  \put(30,30){\bigcircle{60}}
  \put(30,30){\circle*{4}}
  \drawline(6,48)(54,48)
  \drawline(30,0)(20,48)
  \end{picture}
  }
}
\newcommand{\forkTR}
{\raisebox{-28pt}
  {
  \begin{picture}(60,60)
  \put(30,30){\bigcircle{60}}
  \put(30,30){\circle*{4}}
  \drawline(6,48)(54,48)
  \drawline(30,0)(40,48)
  \end{picture}
  }
}
\newcommand{\forkTT}
{\raisebox{-28pt}
  {
  \begin{picture}(60,60)
  \put(30,30){\bigcircle{60}}
  \multiput(12,30)(12,0){4}{\circle*{4}}
  \drawline(12,30)(24,30)
  \drawline(36,30)(48,30)
  \drawline(18,30)(18,2.415)
  \drawline(42,30)(42,2.415)
  \end{picture}
  }
}
\newcommand{\forkDU}
{\raisebox{-28pt}
  {
  \begin{picture}(60,60)
  \put(30,30){\bigcircle{60}}
  \multiput(6,12)(0,36){2}{\line(1,0){48}}
  \qbezier(0,30)(30,30)(30,12)
  \qbezier(60,30)(30,30)(30,48)
  \end{picture}
  }
}
\newcommand{\forkUD}
{\raisebox{-28pt}
  {
  \begin{picture}(60,60)
  \put(30,30){\bigcircle{60}}
  \multiput(6,12)(0,36){2}{\line(1,0){48}}
  \qbezier(0,30)(30,30)(30,48)
  \qbezier(60,30)(30,30)(30,12)
  \end{picture}
  }
}
\newcommand{\forkspoon}
{\raisebox{-28pt}
  {
  \begin{picture}(60,60)
  \put(30,30){\bigcircle{60}}
  \multiput(12,30)(12,0){4}{\circle*{4}}
    \put(12,30){\line(1,0){12}}
    \drawline(18,30)(18,2.415)
    \drawline(48,30)(42,24)
    \qbezier(42,24)(36,18)(24,18)
    \put(24,18){\vector(1,0){0}}
    \qbezier(24,18)(6,18)(6,30)
    \qbezier(6,30)(6,42)(24,42)
    \qbezier(24,42)(36,42)(42,36)
    \drawline(42,36)(48,30)
    \drawline(42,24)(42,2.415)
  \end{picture}
  }
}
\newcommand{\forkstandard}
{
  \begin{picture}(120,120)
  \put(60,60){\bigcircle{120}}
  \multiput(20,60)(16,0){6}{\circle*{4}}
    \multiput(20,60)(32,0){3}{\line(1,0){16}}
    \multiput(28,60)(64,0){2}{\line(0,-1){50.754}}
    \put(60,60){\line(0,-1){60}}
  \end{picture}
}
\newcommand{\forkeight}
{
  \begin{picture}(120,60)
    \multiput(20,40)(80,0){2}{\circle*{4}}
    \put(60,0){\line(0,1){40}}
    \put(20,40){\line(1,0){80}}
    \qbezier(60,40)(75,25)(90,25)
    \put(90,25){\vector(1,0){10}}
    \qbezier(100,25)(120,25)(120,40)
    \qbezier(120,40)(120,60)(100,60)
    \put(100,60){\line(-1,0){20}}
    \qbezier(80,60)(60,60)(50,50)
    \put(50,50){\line(-1,-1){10}}
    \qbezier(40,40)(25,25)(15,25)
    \qbezier(15,25)(0,25)(0,40)
    \qbezier(0,40)(0,60)(20,60)
    \qbezier(20,60)(40,60)(50,50)
    \put(50,50){\line(1,-1){10}}
  \end{picture}
}

\begin{document}

\title{A homological definition of the Jones polynomial}
\author{Stephen Bigelow}
\address{Department of Mathematics,
University of California\\Santa Barbara CA93106, USA}
\email{bigelow@math.ucsb.edu}

\begin{abstract}
We give a new definition of the Jones polynomial.
Let $L$ be an oriented knot or link
obtained as the plat closure of a braid $\beta \in B_{2n}$.
We define a covering space $\tilde{C}$ of
the space of unordered $n$-tuples of distinct points
in the $2n$-punctured disk.
We then describe two $n$-manifolds $\tilde{S}$ and $\tilde{T}$ in $\tilde{C}$,
and show that the Jones polynomial of $L$ can be defined
as an intersection pairing
between $\tilde{S}$ and $\beta\tilde{T}$.
Our construction is similar to one given by Lawrence,
but more concrete.
\end{abstract}

\asciiabstract{We give a new definition of the Jones polynomial.  Let
L be an oriented knot or link obtained as the plat closure of a braid
beta in B_{2n}.  We define a covering space tilde{C} of the space of
unordered n-tuples of distinct points in the 2n-punctured disk.  We
then describe two n-manifolds tilde{S} and tilde{T} in tilde{C}, and
show that the Jones polynomial of L can be defined as an intersection
pairing between tilde{S} and beta tilde{T}.  Our construction is
similar to one given by Lawrence, but more concrete.}


\maketitlepage

\section{Introduction}
\label{sec:intro}

The Jones polynomial is most easily defined using skein relations.
Consider the set of 
formal linear combinations of oriented knots or links in $S^3$ 
over the ring $\Z[q^{\pm \half }]$
modulo the {\em skein relation}
$$q^{-1} L_+ - q L_- = (q^{\half} - q^{-\half}) L_0,$$
where $L_+$, $L_-$ and $L_0$ are oriented knots or links
that are identical except in a ball, where they are as follows.
$$L_+ = \lplus, \; L_- = \lminus, \; L_0 = \lzero.$$
Using this relation,
any oriented knot or link $L$ can be written as 
some scalar multiple of the unknot.
This scalar $V_L \in \Z[q^{\pm \half}]$ 
is uniquely determined by the isotopy class of $L$,
and is called the {\em Jones polynomial} of $L$.

The original definition in \cite{vJ85}
gives the Jones polynomial of the closure of a braid
as a trace function of the image of that braid in the Hecke algebra.
It is natural to ask whether there is a more topological definition -
one which is based not on algebraic properties of a braid,
or combinatorial properties of a projection onto the plane,
but only on topological properties of 
the way the link is embedded into three-dimensional space.

Despite a great deal of effort,
no satisfactory answer to this question is known.
A partial answer was provided by
the groundbreaking paper of Witten \cite{eW89},
which gives an interpretation of the Jones polynomial
in terms of quantum field theory.
Reshetikhin and Turaev \cite{RT91} gave a
mathematically rigorous formulation of this theory,
using quantum groups instead of the Feynman path integral,
which is not yet known to be well-defined.

In this paper,
we follow the approach used by Lawrence in \cite{rL93}.
Let $L$ be the {\em plat closure} of a braid $\beta \in B_{2n}$.
Jones \cite{vJ85} showed that $V_L$ appears as
an entry of the matrix for $\beta$
in a certain irreducible representation of $B_{2n}$.
Lawrence \cite{rL90} gave a topological interpretation of this representation.
This leads to an interpretation of $V_L$
as an intersection pairing between
a certain element of cohomology
and the image under $\beta$ of a certain element of homology.

The definition of $V_L$ given in this paper
is essentially the same as Lawrence's.
However our description of the relevant elements of homology and cohomology
is more explicit.
Indeed we explain how one could use them to directly calculate $V_L$,
which is not clear from \cite{rL90}.
We also give a new and more elementary proof
that our invariant is the Jones polynomial.

This interpretation of $V_L$ might
represent some progress towards a truly topological definition.
However it has a flaw in common with many definitions.
Namely, it is first defined for links in a special form,
and then shown to be an isotopy invariant
by checking it is invariant under certain moves.
In our case, the special form is the plat closure of a braid,
and the moves are those described by Birman in \cite{jB76}.

This paper is an attempt to fill in some details
of a talk given at a workshop in RIMS, Kyoto, in September 2001.
I thank the organisers for their kind hospitality.
This research was supported by the Australian Research Council.

\section{Definitions}

Let $D$ be the unit disk centred at $0$ in the complex plane.
Let $p_1,\dots,p_{2n}$ be points on the real line such that
$$-1 < p_1 < \dots < p_{2n} < 1.$$
Call these {\em puncture points}.  Let
$$D_{2n} = D \setminus \{p_1,\dots,p_{2n}\}.$$
The braid group $B_{2n}$ is the mapping class group of $D_{2n}$.
We will also use other equivalent definitions of $B_{2n}$
as a group of geometric braids,
as the fundamental group of a configuration space of points in the plane,
and as given by generators and relations.
For these and more, see \cite{jB74}.
We use the convention that braids act on $D_{2n}$ on the left,
and geometric braids read from top to bottom.

Let $\bar{C}$ be the set of ordered $n$-tuples of distinct points in $D_{2n}$.
Let $C$ be the quotient of $\bar{C}$ by the symmetric group $S_n$,
that is, the set of unordered $n$-tuples of distinct points in $D_{2n}$.
We now define a homomorphism
$$\Phi \co \pi_1 C \to \langle q \rangle \oplus \langle t \rangle.$$
The motivation for our definition is the fact,
which we will not prove,
that it has a certain universal property.
Namely, any map from $\pi_1 C$ to an abelian group
which is invariant under the action of $B_{2n}$ must factor 
uniquely 
through $\Phi$.

Let $\alpha \co I \to C$ be a loop in $C$.
By ignoring the puncture points we can consider
$\alpha$ as a loop of unordered $n$-tuples in the disk,
and hence as a braid in $B_n$.
Let $b$ be the image of this braid under
the usual abelianisation map from $B_n$ to $\Z$,
which takes each of the the standard generators to $1$.
Similarly, the map
$$s \mapsto \{p_1,\dots,p_{2n}\} \cup \alpha(s)$$
determines a braid in $B_{3n}$.
Let $b'$ be the image of this braid under 
the usual abelianisation map from $B_{3n}$ to $\Z$.
Note that $b$ and $b'$ have the same parity,
equal to the parity of the image of the braid $\alpha$
in the symmetric group $S_n$.
Let $a = \half(b'-b)$.
We define
$$\Phi(\alpha) = q^a t^b.$$

This definition was intended to be easy to state and clearly well-defined,
but it is also somewhat artificial.
A more intuitive definition is as follows.
A loop $\alpha \co I \to C$
can be written as 
$$\alpha(s) = \{\alpha_1(s),\dots,\alpha_n(s)\}$$
for some arcs $\alpha_1,\dots,\alpha_n$ in $D_{2n}$.
The exponent of $q$ in $\Phi(\alpha)$
records the total winding number of these arcs around the puncture points.
The exponent of $t$
records twice the winding number of these arcs around each other.
Thus if two arcs switch places by a counterclockwise half twist
then this contributes a factor of $t$.

Let $\tilde{C}$ be the covering space of $C$ corresponding to $\Phi$.
The group of covering transformations of $\tilde{C}$
is $\langle q \rangle \oplus \langle t \rangle$.
We define the following intersection pairing in $\tilde{C}$.

\begin{defn}
Suppose $A$ and $B$ are immersed $n$-manifolds in $C$
such that at least one of them is closed,
and the other intersects every compact subset of $C$
in a compact set.
Suppose $\tilde{A}$ and $\tilde{B}$ are lifts of
$A$ and $B$ respectively to $\tilde{C}$.
For any $a,b\in\Z[q^{\pm 1},t^{\pm 1}]$ let $q^a t^b \tilde{A}$ be 
the image of $\tilde{A}$ under the covering transformation $q^a t^b$.
There is a well-defined algebraic intersection number 
$(q^a t^b \tilde{A}, \tilde{B}) \in \Z$.
We define
$$\langle \tilde{A},\tilde{B} \rangle = \sum_{a,b\in\Z}
(q^a t^b\tilde{A},\tilde{B})q^a t^b \in \Z[q^{\pm 1},t^{\pm 1}].$$
This is a finite sum,
since the number of nonzero terms is at most
the geometric intersection number of $A$ and $B$ in $C$.
\end{defn}

To specify an $n$-manifold $\tilde{C}$,
it will help to have a fixed basepoint.

\begin{defn}
Let $d_1,\dots,d_n$ be distinct points on $\partial D$.
We take them to lie in the lower half plane, ordered from left to right.
Let ${\bf c} = \{d_1,\dots,d_n\}$
and fix a choice of $\tilde{\bf c} \in \tilde{C}$ in the fibre over ${\bf c}$.
\end{defn}

We now define a certain type of picture in the disk
which we will use to represent an immersed $n$-manifold in $\tilde{C}$.

\begin{defn}
A {\em fork diagram} in $D_{2n}$ consists of maps
$$E_1,\dots,E_n \co I \to D$$
called {\em tine edges}, and maps
$$E'_1,\dots,E'_n \co I \to D$$
called {\em handles}, subject to the following conditions.
\begin{itemize}
\item the tine edges are disjoint embeddings
  of the interior of $I$ into $D_{2n}$,
  and map the endpoints of $I$ to the puncture points in $D$
  (not necessarily injectively),
\item the handles are disjoint embeddings of $I$ into $D_{2n}$,
\item $E'_i$ is a path from $d_i$
  to a point in the interior of $E_i$.
\end{itemize}
\end{defn}

Such a fork diagram determines an immersed open $n$-ball 
$\tilde{U}$ in $\tilde{C}$ as follows.
Note that $E_1 \times \dots \times E_n$
maps the interior of $I \times \dots \times I$ into $\bar{C}$.
Now let $U$ be the projection of this map to $C$.
Let $\gamma$ be the path in $C$ given by
$$\gamma(s) = \{E'_1(s),\dots,E'_n(s)\}.$$
Lift this to a path 
$\tilde{\gamma}$ in $\tilde{C}$ starting at $\tilde{\bf c}$.
We define $\tilde{U}$ to be
the lift of $U$ which contains $\tilde{\gamma}(1)$.
This is an oriented open $n$-ball in $\tilde{C}$.

\begin{figure}
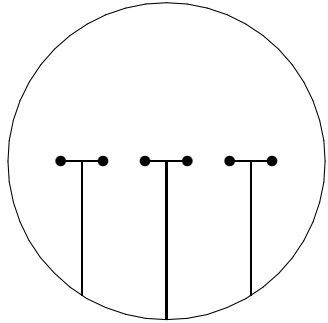

\centering
\forkstandard
\caption{The standard fork diagram for $n=3$}
\label{fig:standard}
\end{figure}
\begin{defn}
Let $\tilde{S}$ denote the open $n$-ball in $\tilde{C}$
corresponding to the {\em standard fork diagram}
shown in Figure \ref{fig:standard},
where all tine edges are oriented from left to right.
\end{defn}

\begin{figure}
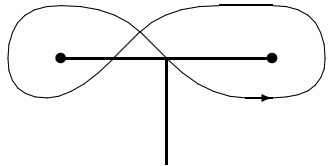

\centering
\forkeight
\caption{The figure eight $F_i$ corresponding to $E_i$}
\label{fig:eight}
\end{figure}
\begin{defn}
For each tine edge $E_i$ in the standard fork diagram,
let $F_i$ be the map from $S^1$ to the figure-eight
as shown in Figure \ref{fig:eight}.
Then $F_1 \times \dots \times F_n$
is an immersion from the $n$-torus into $\bar{C}$.
Let $T$ be the projection of this map into $C$.
Use the handles to specify a lift $\tilde{T}$ of $T$ to $\tilde{C}$,
as in the definition of $\tilde{S}$.
\end{defn}

If $[\beta] \in B_{2n}$,
let $\beta$ be a homeomorphism from $D_{2n}$ to itself
that represents the mapping class $[\beta]$.
Let $\beta'$ be the induced map from $C$ to itself.
This can be lifted to a map from $\tilde{C}$ to itself.
Let $\tilde{\beta}'$ be the lift of $\beta'$ which fixes $\tilde{\bf c}$.
By abuse of notation
we will use $\beta$ to denote $[\beta]$, $\beta'$ and $\tilde{\beta}'$.

\begin{defn}
An {\em oriented braid} is a braid $\beta \in B_{2n}$
together with a choice of orientation for each of the $2n$ strands
such that the orientations match up correctly
when we take the plat closure of $\beta$.
\end{defn}

We now define our invariant $V'_\beta$
of an oriented braid $\beta \in B_{2n}$.

\begin{defn}
Suppose $\beta \in B_{2n}$ is an oriented braid.
Let $w$ be the writhe of $\beta$,
that is, the number of right handed crossings 
minus the number of left handed crossings.
Let $e$ be the sum of the exponents of the generators
in any word representing $\beta$.
Let
$$\lambda = (-q^{\half} - q^{-\half})^{-1}
            (-q^{\frac{3}{4}})^w (q^{-\frac{1}{4}})^{e+2n}.
$$
Then we define
$$V'_\beta = \lambda \langle \tilde{S}, \beta\tilde{T} \rangle
           |_{(t=-q^{-1})}.$$
\end{defn}

The main result of this paper is the following.

\begin{thm}
\label{thm:main}
$V'_\beta$ is the Jones polynomial of the plat closure of $\beta$.
\end{thm}

Note that the simpler formula
$$q^{-\frac{1}{4}(e+2n)} \langle \tilde{S},\beta\tilde{T} \rangle
           |_{(t=-q^{-1})}$$
therefore gives the Kauffman bracket of the plat closure of $\beta$,
normalised to equal one for the empty diagram.

\section{How to compute the invariant}
\label{sec:compute}

In this section, we describe of how one could compute 
$\langle \tilde{S},\beta\tilde{T} \rangle$, and hence $V'_\beta$.

Let $E_1,\dots,E_n$ and $E'_1,\dots,E'_n$
be the tine edges and handles of the standard fork diagram.
Let $F_1,\dots,F_n$ be the corresponding figure-eights
as in Figure \ref{fig:eight}.
By applying an isotopy,
we can assume that each $\beta F_i$ and $E_j$ intersect transversely
and there are no triple points.

The intersection $S \cap \beta T$ consists of points
${\bf e} = \{e_1,\dots,e_n\}$ in $C$
such that $e_i \in E_i \cap F_{\pi i}$
for some permutation $\pi$ of $\{1,\dots,n\}$.
Each such ${\bf e}$ contributes a monomial
$\pm q^a t^b$ to $\langle \tilde{S},\beta\tilde{T} \rangle$.
The sign of this monomial is the sign of the intersection
of $S$ with $\beta T$ at ${\bf e}$.
Let $\tilde{\bf e}$ be the lift of ${\bf e}$ which lies in $\tilde{S}$.
The integers $a$ and $b$ are such that 
$q^a t^b \tilde{\bf e}$ lies in $\beta\tilde{T}$.
The sum of these terms $\pm q^a t^b$ over all ${\bf e}$
in $S \cap \beta T$
is the required polynomial $\langle \tilde{S},\beta\tilde{T} \rangle$.

We now describe how to compute the monomial for a given ${\bf e}$ as above.
Let $m$ be the number of points $e_i$
such that the sign of the intersection
of $E_i$ with $F_{\pi i}$ at $e_i$ is negative.
The sign of our required monomial is then
$(-1)^m$ times the parity of the permutation $\pi$.

Let $h$ be the path
$$h(s) = \{E'_1(s),\dots,E'_n(s)\}$$
in $C$.
Let $\epsilon$ be the path
$$\epsilon(s) = \{\epsilon_1(s),\dots,\epsilon_n(s)\},$$
where $\epsilon_i$ is a segment of $E_i$ going from $E'_i(1)$ to $e_i$.
Then the composition $h\epsilon$ 
lifts to a path from $\tilde{\bf c}$ to $\tilde{\bf e}$.
Let $\phi$ be the path
$$\phi(s) = \{\phi_1(s),\dots,\phi_n(s)\},$$
where $\phi_i$ is a segment of $\beta F_{\pi i}$
going from $\beta E'_{\pi i}(1)$ to $e_i$.
Then the composition $\beta(h) \phi$
lifts to a path from $\tilde{\bf c}$ to $q^a t^b \tilde{\bf e}$.
Recall that $q^a t^b \tilde{\bf e} \in \beta\tilde{T}$.
We conclude that the path
$$\delta = \beta(h) \phi \epsilon^{-1} h^{-1}$$
lifts to a path from $\tilde{\bf c}$ to $q^a t^b \tilde{\bf c}$.
Thus
$$q^a t^b = \Phi(\delta).$$

It is possible to calculate the Jones polynomial
of a knot or link by these methods.
Given a knot diagram,
replace every sequence of consecutive undercrossings with a figure-eight.
Attach handles in any convenient fashion
and proceed as described above.
Then either compute the correct factor $\lambda$,
or simply settle for the value of the Jones polynomial
up to sign and multiplication by a power of $q^{\half}$.
I have performed this calculation by hand
for the unknot, the Hopf link, and the trefoil knot.
Even for the trefoil, it was necessary to correct several mistakes
before reaching the correct answer.
A computer implementation would be more reliable,
but probably no more efficient than existing methods.

\section{Lemmas}

The pairing $\langle \cdot,\cdot \rangle$
is really a pairing between $H^n(\tilde{C})$ and $H_n(\tilde{C})$.
These cohomology and homology groups
are modules over $\Z[q^{\pm 1},t^{\pm 1}]$,
where $q$ and $t$ act by covering transformations.
The open $n$-ball corresponding to a fork diagram
represents an element of $H^n(\tilde{C})$.
We can thus extend the algebraic pairing
to take linear combinations of fork diagrams in its first entry.
We now prove some relations that hold between fork diagrams
considered as elements of $H^n{\tilde{C}}$.

\begin{lem}
\label{lem:relations}
The following relations hold.
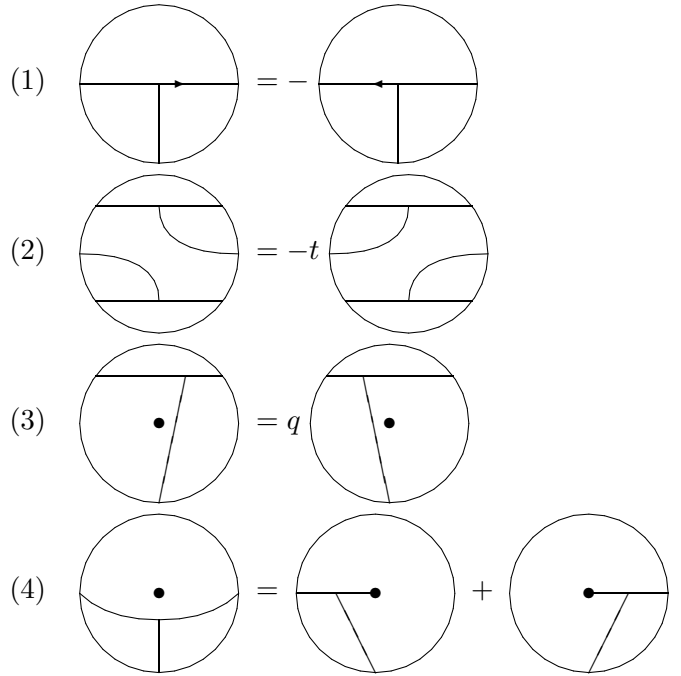
\begin{enumerate}
\item $\forkTplus = -\forkTminus$
\item $\forkDU = -t\forkUD$
\item $\forkTR = q\forkTL$
\item $\forkY = \forkL + \forkR$.
\end{enumerate}
\end{lem}

Here, the diagrams in each relation are understood to be identical
except in the disk shown.
The disks shown are also allowed to contain
an arbitrary number of additional tine edges,
which are required to be identical for all diagrams in the relation.
All tine edges are oriented in any consistent fashion
except in relation (1).

\begin{proof}[Proof of Lemma \ref{lem:relations}]
In the first three relations,
the tine edges are the same in both sides of the equation.
Thus the corresponding open $n$-balls are lifts of 
the same open $n$-ball in $C$,
possibly with different orientations.

To determine the orientation,
recall that the open $n$-ball in $C$ was defined
as the product $T_1\times \dots \times T_n$ of tine edges,
where $T_i$ is the tine edge whose handle is attached to $d_i$.
Thus the orientation is determined by
the orientations of the tine edges
and the order in which they occur.
In the first relation, the orientation of one tine edge was reversed.
In the second, the order in which they occur underwent a transposition.
In the third, there was no change to
the order or orientation of the tine edges.
Thus the signs are as claimed in these relations.

The choice of lift is determined by the handles.
In the first relation, the handles are unchanged.
In the second and third, the changes in the handles
correspond to the scalars $t$ and $q$ as given.

The final relation describes the process of
cutting the open $n$-ball along a hyperplane
by pushing that hyperplane into a puncture point.
\end{proof}

One useful consequence of these relations is the following relation.
(Thanks to Saul Schleimer for suggesting the name.)

\begin{lem}[The fork-spoon relation]
\label{lem:forkspoon}
Modulo $(1+qt)$, we have
$$\forkspoon = (1-q^{-1})\forkTT,$$
where the horizontal tine edges are oriented left to right.
\end{lem}

\section{Proof of the main theorem}

Throughout this section we use the conventions that $t=-q^{-1}$,
and the horizontal tine edges in any fork diagram
are oriented from left to right.

Let $\beta$ be an oriented braid
and let $L$ be the plat closure of $\beta$.
The aim of this section is to prove Theorem \ref{thm:main},
that $V'_\beta = V_L$.
First we prove that
$V'_\beta$ depends only on the isotopy class of $L$.
We use a result due to Birman
which gives the ``Markov moves'' for plat closures.

\begin{defn}
Let $K_{2n}$ be the subgroup of $B_{2n}$ generated by
\begin{itemize}
\item $\sigma_1$,
\item $\sigma_2 \sigma_1^2 \sigma_2$,
\item $\sigma_{2i} \sigma_{2i-1} \sigma_{2i+1} \sigma_{2i}$
      for $i=1,\dots,n-1$,
\end{itemize}
where $\sigma_1,\dots,\sigma_{2n-1}$ are 
the standard generators of $B_{2n}$.
\end{defn}

\begin{lem}[Birman]
\label{lem:birman}
Two oriented braids 
$\beta_1 \in B_{2n_1}$ and $\beta_2 \in B_{2n_2}$
have isotopic plat closures if and only if 
they are related by a finite sequence of the following moves.
\begin{itemize}
\item $\beta \mapsto g \beta h$,
      where $\beta \in B_{2n}$ and $g,h \in K_{2n}$.
\item $\beta \leftrightarrow \sigma_{2n}\beta$,
      where $\beta \in B_{2n}$ and $\sigma_{2n}\beta \in B_{2n+2}$.
\end{itemize}
\end{lem}

The move $\beta \mapsto \sigma_{2n}\beta \in B_{2n+2}$ is called
{\em stabilisation}.

Our statement of Lemma \ref{lem:birman}
differs from that of \cite{jB76} in two respects.
First, Birman's result applies only to knots.
This was necessary in order to show that 
the two types of move commute.
We do not need the moves to commute,
so we can apply Birman's proof to the case of links.

Second, Birman did not consider the issue of orientation.
There is a small technical point here.
Our statement of the lemma should really specify
the effect of each move on the orientation of the braid.
However since each move corresponds to an isotopy of the plat closure,
it is clear what this effect should be.
With this in mind, Birman's proof goes through unchanged.

\begin{lem}
\label{lem:Kbeta}
If $\beta \in B_{2n}$ is an oriented braid and $g \in K_{2n}$
then $V'_{g \beta} = V'_\beta$.
\end{lem}

\begin{proof}
It suffices to prove this
in the case $g$ is one of the generators of $K_{2n}$.

Recall that the coefficient in the definition of $V'_\beta$ is
$$\lambda = (-q^{\half} - q^{-\half})^{-1}
            (-q^{\frac{3}{4}})^w (q^{-\frac{1}{4}})^{e+2n},
$$
where $w$ and $e$ are respectively the writhe and exponent sum of $\beta$.
Now $\sigma_1$ contributes a left-handed crossing to $\sigma_1\beta$,
so the writhe of $\sigma_1 \beta$ is $w-1$.
For all other generators, $g\beta$ has writhe $w$.
Thus we must show that
$$\delta \langle \tilde{S}, g\beta\tilde{T} \rangle
  = \langle \tilde{S}, \beta\tilde{T} \rangle,$$
where $\delta$ is $-q^{-1}$ for $g=\sigma_1$,
and $q^{-1}$ for all other generators.
We can rewrite this as
$$\langle \delta^{-1}\tilde{S}, g\beta\tilde{T} \rangle
  = \langle g\tilde{S}, g\beta\tilde{T} \rangle.$$
Thus it suffices to prove the identities
$$g\tilde{S} = \delta^{-1} \tilde{S}.$$
In the cases $g=\sigma_1$ and
$g=\sigma_{2i}\sigma_{2i-1}\sigma_{2i+1}\sigma_{2i}$,
this follows easily from Lemma \ref{lem:relations}.
In the case $g=\sigma_2\sigma_1^2\sigma_2$,
it helps to first use the fork-spoon relation.
\end{proof}

\begin{lem}
\label{lem:betaK}
If $\beta \in B_{2n}$ is an oriented braid and $h \in K_{2n}$
then $V'_{\beta h} = V'_\beta$.
\end{lem}

\begin{proof}
Note that $\tilde{T}$ can be isotoped to be equal to
$(1-q)^n\tilde{S}$ except in
an arbitrarily small neighbourhood of the puncture points.
Thus
\begin{eqnarray*}
\langle \tilde{S},\beta\tilde{T} \rangle 
&=& (1-q)^{-n} \langle \tilde{T}, \beta\tilde{T} \rangle \\
&=& (1-q)^{-n} \langle \beta^{-1}\tilde{T}, \tilde{T} \rangle \\
&=& \langle \beta^{-1}\tilde{T},\tilde{S} \rangle.
\end{eqnarray*}
Similarly
$$\langle \tilde{S},\beta h\tilde{T} \rangle
 = \langle \beta^{-1}\tilde{T}, h\tilde{S} \rangle.
$$
The result now follows from the proof of the previous lemma.
\end{proof}

Before we move on to stabilisation,
we prove the following special case of the skein relation.

\begin{lem}
\label{lem:TRskein}
Let $\beta\in B_{2n}$ be an oriented braid.
Suppose the strands second and third from the right at the top
have a parallel orientation.
Let $\beta_+ = \sigma_{2n-2}\beta$ and $\beta_- = \sigma_{2n-2}^{-1}\beta$.
Then
$$q^{-1}V'_{\beta_+} - qV'_{\beta_-}
  = (q^{\half} - q^{-\half})V'_\beta.$$
\end{lem}

\begin{proof}
If the writhe of $\beta$ is $w$
then the writhe of $\beta_+$ is $w+1$
and the writhe of $\beta_-$ is $w-1$.
A simple calculation reduces the problem to showing that
$$(\sigma_{2n-2}-1)(\sigma_{2n-2}+q)\tilde{S} = 0.$$
By the fork-spoon lemma, this is equivalent to
$$(\sigma_{2n-2}-1)(\sigma_{2n-2}+q)\forkspoon = 0,$$
where we have shown only a disk containing
the puncture points $p_{2n-3},\dots,p_{2n}$.
To simplify notation,
assume that $n=2$
and show only a disk containing $p_1$, $p_2$ and $p_3$.
This reduces the problem to showing that
$$(\sigma_2-1)(\sigma_2+q)\forkTo = 0.$$
By Lemma \ref{lem:relations},
$$(\sigma_2-1)\forkTo = \forkoT,$$
and
$$(\sigma_2+q)\forkoT = 0.$$
This completes the proof.
\end{proof}

\begin{lem}
Suppose $\beta \in B_{2n}$ is an oriented braid.
Let $\beta'_+ = \sigma_{2n}\beta \in B_{2n+2}$.
Then $V'_{\beta'_+} = V'_\beta$.
\end{lem}

\begin{proof}
Let $\beta'$ be the image of $\beta$ in $B_{2n+2}$.
Note that $\beta'_+ = \sigma_{2n}\beta'$.
Let $\beta'_- = \sigma_{2n}^{-1}\beta'$.
Then 
$$\beta'_+ \sigma_{2n+1}^2 = \sigma_{2n}\sigma_{2n+1}^2\sigma_{2n} \beta'_-.$$
One can show that $\sigma_{2n+1}$ and $\sigma_{2n}\sigma_{2n+1}^2\sigma_{2n}$
lie in $K_{2n+1}$.
By Lemmas \ref{lem:Kbeta} and \ref{lem:betaK}, it follows that
$$V'_{\beta'_+} = V'_{\beta'_-}.$$
By Lemma \ref{lem:TRskein},
$$q^{-1}V'_{\beta'_+} - qV_{\beta'_-}
  = (q^{\half} - q^{-\half})V'_{\beta'}.$$
Thus
$$(-q^{\half}-q^{-\half})V'_{\beta'_+} = V'_{\beta'}.$$
To show that $V'_{\beta'_+} = V'_\beta$,
it therefore suffices to show that 
$$V'_{\beta'} = (-q^{\half}-q^{-\half})V'_\beta.$$
Note that $\beta$ and $\beta'$ have the same writhe $w$,
and the same exponent sum $e$,
but have $2n$ and $2n+2$ strands respectively.
Thus coefficient $\lambda$ used in the computation of $V'_{\beta'}$ 
is $q^{-\half}$ times that of $V'_\beta$.
The problem is therefore reduced to showing that
$$\langle \tilde{S'},\beta\tilde{T'} \rangle
  = (-1-q) \langle \tilde{S},\beta\tilde{T} \rangle,$$
where $\tilde{S'}$ is the open $(n+1)$-ball
corresponding to the standard fork diagram in $D_{2n+2}$,
and $\tilde{T'}$ is the corresponding $(n+1)$-torus.

Now $\tilde{S'}$ is the product of $\tilde{S}$ with an edge,
and $\tilde{T'}$ is the product of $\tilde{T}$ with a circle.
This edge meets this circle at two points in the disk.
These two points contribute
$-1$ and $-q$ times $\langle \tilde{S},\beta\tilde{T} \rangle$.
To see this, consider the computation of the pairing
described in Section \ref{sec:compute}.
\end{proof}

We have shown that $V'_\beta$ is invariant
under each of the moves in \ref{lem:birman}.
Thus it is an isotopy invariant of the plat closure $L$ of $\beta$.

We now prove the skein relation given in Section \ref{sec:intro}.
Let $L_+$, $L_-$ and $L_0$ be as defined there.
Since $V'$ is an isotopy invariant,
we are free to choose any representation of these links
as plat closures of a braid.
In particular, we can move the ball on which the three links differ
to the top right of the diagram.
Then we can isotope the rest of the diagram to be a plat closure.
Thus we can reduce to the case which was already proved
in Lemma \ref{lem:TRskein}.

Finally,
a direct computation as described in Section \ref{sec:compute}
verifies that our invariant of the unknot is one.
This completes the proof that $V'_\beta$ is the Jones polynomial
of the plat closure of $\beta$.


\Addresses\recd

\end{document}